\documentclass[a4paper,11pt]{article}

\usepackage{latexsym}
\usepackage{amssymb}
\usepackage{amsmath}
\usepackage[dvips]{graphicx}

\title{\LARGE\bf Correction of Errors in the First-Order Perturbation Expansions of Singular Vectors}

\author {Gabriel Ok\v sa \\
        Slovak Academy of Sciences, Institute of Mathematics\\
        Bratislava, Slovak Republic}

\newcommand{\old}[1]{{}}

\begin{document}

\date{}

\maketitle

This note provides a correction of several serious errors (not only misprints!) 
in \cite[p.208]{stew}, where Stewart provides incorrect equations for the first-order expansions 
of left and right singular vectors of a perturbed $n\times p$ matrix $\tilde{X} = X+E$. 
Stewart's notation is used throughout.

In \cite[(3.6), p.207]{stew}, only two elements of the perturbed matrix on the right-hand side are given 
explicitly. Here, all of them are listed as follows:
\begin{align} \label{pert}  
  \varphi_1 &= u_1^T E v_1,\quad f_{12}^T = u_1^T E V_2,\quad f_{21} = U_2^T E v_1, \nonumber \\
  f_{31} &= U_3^T E v_1,\quad F_{22} = U_2^T E V_2,\quad F_{32} = U_3^T E V_2.   
\end{align}       
First two errors are in expressions for $g_2$ and $h_2$ in \cite[(3.11), p.208]{stew}. These vectors are solutions 
of the system of equations \cite[(3.9)--(3.10),p.208]{stew}, but they are wrong. To understand why, take \cite[(3.9), p.208]{stew},  
and write
\begin{equation*}
   g_2 = \sigma_1^{-1} (f_{21} + \Sigma_2 h_2).
\end{equation*}
Plug it into (3.10),
\begin{equation*}
  \sigma_1 h_2 - \Sigma_2 \left (  \sigma_1^{-1} (f_{21} + \Sigma_2 h_2) \right ) = f_{12}
\end{equation*}
and compute
\begin{equation}   \label{h2}
  h_2 = (\sigma_1^2 I - \Sigma_2^2)^{-1}\, (\sigma_1\, f_{12} + \Sigma_2\, f_{21}).
\end{equation}
Then the back substitution gives
\begin{multline*}
  g_2 = \sigma_1^{-1} f_{21} + (\sigma_1^2 I - \Sigma_2^2)^{-1} \Sigma_2 f_{12} + 
            (\sigma_1^2 I - \Sigma_2^2)^{-1} \Sigma_2^2 \sigma_1^{-1} f_{21} \\
        = (\sigma_1^2 I - \Sigma_2^2)^{-1} \Sigma_2 f_{12} + [I + (\sigma_1^2 I - \Sigma_2^2)^{-1} \Sigma_2^2] \sigma_1^{-1} f_{21}.
\end{multline*}
But
\begin{align*}
   I + (\sigma_1^2 I - \Sigma_2^2)^{-1} \Sigma_2^2 &= \mathrm{diag} \left (  1+ \frac{\sigma_k^2}{\sigma_1^2 - \sigma_k^2} \right ) 
        = \mathrm{diag} \left (  \frac{\sigma_1^2}{\sigma_1^2 - \sigma_k^2} \right ) \\
        &= \sigma_1^2 (\sigma_1^2 I - \Sigma_2^2)^{-1},
\end{align*}
so that 
\begin{equation}   \label{g2}
   g_2 = (\sigma_1^2 I - \Sigma_2^2)^{-1}\, (\sigma_1\, f_{21} + \Sigma_2\, f_{12}).
\end{equation}
Hence, the correct solution presented in Eqs.~(\ref{h2})--(\ref{g2})  differs from that in \cite[(3.11),p.208]{stew} 
by plus signs (instead of minus signs)  
after $f_{12}$ and $f_{21}$.
These errors contaminate also Stewart's Theorem 3.4 on p.208 in \cite{stew}. But this theorem contains three 
additional errors. To understand them, the detailed derivation of equations 
for $\tilde{u}_1$ and $\tilde{v}_1$ is provided now; this derivation is not given in the book.

Write \cite[(3.7), p.207]{stew}, in modified form:
\begin{equation*}
  U^T \tilde{X}V \begin{pmatrix} 1 \\ h_2  \end{pmatrix} = (\sigma_1 + \theta_1)
                 \begin{pmatrix} 1 \\ g_2 \\ g_3   \end{pmatrix},
\end{equation*}  
which is equivalent to
\begin{equation*}
    \tilde{X} V \begin{pmatrix} 1 \\ h_2  \end{pmatrix} = (\sigma_1 + \theta_1) 
                 U \begin{pmatrix} 1 \\ g_2 \\ g_3   \end{pmatrix}. 
\end{equation*}
This means that
\begin{equation*}
   \tilde{u}_1 = U \begin{pmatrix} 1 \\ g_2 \\ g_3   \end{pmatrix}, \quad 
   \tilde{v}_1 = V \begin{pmatrix} 1 \\ h_2  \end{pmatrix}.
\end{equation*}
Using partitions $U = (u_1, U_2, U_3),\  V=(v_1, V_2),$ (see \cite[p.206]{stew}, bottom) together 
with Eq.~(\ref{pert}) leads to 
\begin{align}  
  \begin{split}
  \tilde{u}_1 &= u_1 + U_2 g_2 + U_3 g_3 + O(\|E\|^2) \\  \nonumber
     &= u_1 + U_2 (\sigma_1^2 I - \Sigma_2^2)^{-1} (\sigma_1U_2^T E v_1 + \Sigma_2 V_2^T E^T u_1) \\ 
      &\quad  + \sigma_1^{-1} U_3\, U_3^T\, E\, v_1  + O(\|E\|^2), 
  \end{split}    \\
  \begin{split}  \label{cor}
   \tilde{v}_1 &= v_1 + V_2 h_2 + O(\|E\|^2) \\ 
    &= v_1 + V_2 (\sigma_1^2 I - \Sigma_2^2)^{-1} (\sigma_1V_2^T E^T u_1 + \Sigma_2 U_2^T E v_1) \\
    &\quad   + O(\|E\|^2).   
 \end{split}
\end{align}
In summary, a comparison of Eq.~(\ref{cor}) with \cite[Th.3.4,(3.12),p.208]{stew} reveals the following errors in 
the book:
\begin{enumerate}
 \item In the formula for $\tilde{u}_1$:
   \begin{itemize} 
     \item Wrong sign after $v_1$ (minus instead of plus).
 
     \item Application of $V_2$ instead of $V_2^T$ after $\Sigma_2$.
     
     \item Total omission of $U_3$ after $\sigma_1^{-1}$.
   \end{itemize}  
  \item In the formula for $\tilde{v}_1$:
    \begin{itemize}
      \item Application of $V_2$ instead of $V_2^T$ after $\sigma_1$.
    
      \item Wrong sign after $u_1$ (minus instead of plus).             
    \end{itemize}
\end{enumerate}
%
\section*{Acknowledgenent}
This research was supported by the VEGA grant no.~2/7143/27.


\end{document}